\documentclass[12pt]{amsart}
\input epsf
\usepackage{amsmath}
\usepackage{amssymb}
\usepackage{amscd}
\usepackage{amsthm}
\usepackage{subfigure}

\usepackage[all]{xy}
\usepackage{color}

\numberwithin{equation}{section}

\newtheorem{theorem}[equation]{Theorem}

\newtheorem{proposition}[equation]{Proposition}

\newtheorem{lem}[equation]{Lemma}
\newtheorem{lemma}[equation]{Lemma}

\newtheorem{corollary}[equation]{Corollary}

\theoremstyle{remark}

\theoremstyle{definition}
\newtheorem{definition}[equation]{Definition}

\def\Xint#1{\mathchoice 
	{\XXint\displaystyle\textstyle{#1}}%
	{\XXint\textstyle\scriptstyle{#1}}%
	{\XXint\scriptstyle\scriptscriptstyle{#1}}%
	{\XXint\scriptscriptstyle\scriptscriptstyle{#1}}%
	\!\int} 
\def\XXint#1#2#3{{\setbox0=\hbox{$#1{#2#3}{\int}$} 
	\vcenter{\hbox{$#2#3$}}\kern-.5\wd0}} 
 
\def\dashint{\Xint-}

\newcommand{\ra}{\rightarrow}

\newcommand{\N}{\mathbb N}
\newcommand{\n}{{\mathcal N}}

\newcommand{\R}{\mathbb R}

\renewcommand{\H}{\mathcal H}
\renewcommand{\L}{{\mathcal L}}

\newcommand{\F}{{\mathcal F}}

\newcommand{\diam}{\operatorname{diam}}

\newcommand{\bilip}{biLipschitz }

\newcommand{\Length}{\operatorname{Length}}

\newcommand{\rips}{\operatorname{Rips}}

\newcommand{\sur}{\operatorname{Sur}}

\def\de{\delta}

\def\eps{\epsilon}
\def\ga{\gamma}

\def\ra{\rightarrow}

\def\si{\sigma}
\def\Si{\Sigma}

\def\defeq{:=}

\begin{document}

\title{Doubling property for 
biLipschitz homogeneous geodesic surfaces}

\author{Enrico Le Donne}

\date{March 1, 2010}

\begin{abstract}In this paper we discuss general properties of geodesic surfaces that are locally biLipschitz homogeneous. In particular, we prove that they are locally doubling and that there exists a special doubling measure analogous to the Haar measure for locally compact groups.
\end{abstract}

\maketitle

\section{Introduction}
According to a consequence of a general theorem by V. N. Berestovski\u\i\; \cite{b, b1, b2}, if a geodesic distance $d$ on a surface $S$ induces the surface topology of $S$ and has the property that the isometries of $(S,d)$ act transitively on $S$, 
then $(S,d)$ is isometric to a Finsler surface. In particular, such spaces are locally \bilip equivalent to a planar Euclidean domain.

Although, some geodesic distances on the plane are not locally \bilip equivalent to the Euclidean distance. Laakso constructed in \cite{Laakso} geodesic metrics on the plane that are not \bilip embeddable into any $\R^n$, but still share many properties with the Euclidean metric. Some of these properties are   Ahlfors $2$-regularity, local linear contractibility, and the fact that a Poincar\'e inequality holds; see \cite{Heinonenbook} for an introduction to these last definitions.

In this paper we begin the study of a property that holds in the case of the Euclidean plane but has never been singled out: the fact that \bilip maps act transitively.
Since 
every Riemannian/Finsler surface is locally \bilip equivalent to an Euclidean planar domain, 
every two points on the surface have neighborhoods that are \bilip equivalent. Briefly, we say that every Finsler surface is locally \bilip homogeneous; see the next section for the general definitions.
Thus our natural question   is whether every geodesic distance on the plane, or on a surface, where the \bilip maps act locally transitively, is \bilip equivalent to a Riemannian distance and so, locally, to the Euclidean distance.

General homogeneity appears frequently in different mathematical areas and is as natural to assume as it is hard to handle in proofs. We refer, for example, to the challenging open conjecture of Bing and Borsuk, \cite{bing-borsuk}, which states that an $n$-dimensional, homogeneous,  absolute neighborhood retract,  should be an $n$-manifold. See  \cite{Bryant, Halverson-Repovs} for definitions, progress and references.

Homogeneity by isometries in the case of geodesic metric spaces has been successfully studied and characterized by Berestovski\u\i\; \cite{b, b1, b2}. The interest in \bilip homogeneity is  relatively recent. It has been studied by several authors  \cite{Bishop, GH, Freeman-Herron} in dimension one   for planar curves with metrics induced by the ambient geometry. BiLipschitz homogeneity for geodesic spaces has appeared naturally in Geometric Group Theory for some actions on quasi-planes, i.e., geometric objects that are coarsely $2$ dimensional, e.g., in \cite{Kapovich-Kleiner, Kapovich-Kleiner2}.

Our purpose is to study the $2$-dimensional case together with the hypothesis, as is common  in Geometric Group Theory, that the metric is geodesic. Such an assumption in dimension one would give trivial results.

The main result of this paper is  that any geodesic metric surface that is locally \bilip homogeneous is a locally doubling metric space. 
This fact leads to plenty of consequences, e.g., the Hausdorff dimension is finite and there exists a doubling measure that, like the Haar measure on Lie groups is preserved by (left) translations, 
 is ``\bilip preserved'' by \bilip maps.
 
\subsection{Definitions, results, and strategies}
In a metric space $(X,d)$, the {\it length of a curve} $\gamma : [a, b] \rightarrow X$ is
$$
     \Length_d(\gamma):=\sup \left\{ \sum_{i=1}^n d(\gamma(t_i),\gamma(t_{i-1})) : n \in \mathbb{N} \mbox{ and } a = t_0 < t_1 < \cdots < t_n = b \right\}. 
    $$ 
 A {\it rectifiable curve} is a curve with finite length. 
A {\em geodesic  space} is a metric space where any two points are the end points of a rectifiable curve whose length is exactly the distance between the two points.

A metric space $(X,d)$ is {\em doubling} if
there is a constant $N\in \N$ such that each ball
$B(x,2R)\subset X$  is contained in
the union of $\le N$ balls of radius $R$.
We say that $(X,d)$ is {\em locally doubling} if any point has a neighborhood that is doubling.

We say that a  metric space $(X,d)$ is {\em locally  biLipschitz  homogeneous} 
if, for every two points
$x_1,x_2 \in X$, there are neighborhoods $U_1$ and $U_2$  of $x_1$ and $x_2$ respectively
and a \bilip homeomorphism $f:U_1\to U_2$, such that $f(x_1)=x_2$.

A metric space $(X,d)$ is   {\em locally linearly contractible} if there is a
constant $C\geq1$ such that each metric ball of radius $ R < C^{ -1}$
in the space can be 
contracted to a point inside the ball of the same center but radius $C R$. See \cite{Semmes96} for an ample 
analysis of this condition. 


We prove the following: 
\begin{theorem}\label{teorema1}
Let $(X,d)$ be a geodesic metric space topologically equivalent to a surface.
Assume that $X$ is 
 locally {\bilip  homogeneous}.
Then \begin{enumerate}
\item The metric space $(X,d)$ is  locally doubling,
\item The Hausdorff dimension of  $(X,d)$ is finite,
\item The Hausdorff $2$-measure $\H^2$ of small $r$-ball $B_r$ has a quadratic lower bound: for each point $p\in X$, there are constants $c,\bar r>0$, so that
$$\qquad \H^2(B(p,r))\geq c r^2, \quad \text{ for } r<\bar r,$$  and
\item Every point of $(X,d)$ has a neighborhood that is  locally linearly contractible.
\end{enumerate} \end{theorem}

Subsequently, we investigate the properties of general doubling \bilip homogeneous spaces.
We show that they admit an analog of the Haar measure: 
there exists a doubling measure that is quasi-preserved    
by \bilip maps and is quasi-unique. 
For $\alpha>0$, we consider two Borel measures
  $\nu$ and $\mu$ to be $\alpha$-quasi-equivalent,
writing $\nu\stackrel{\alpha}{\approx} \mu$,
when, for all Borel sets $A$,
  \begin{eqnarray}
 \frac{1}{\alpha} \mu(A)\le \nu(A))\le\alpha  \mu(A).
  \end{eqnarray}
  With such a  notation we can precisely formulate the result.
 \begin{proposition} [] [Existence]\label{Haar measure} Let $X$ be a locally compact and separable metric space whose metric is  doubling. 
Then there exists a (non-zero) Radon measure $\mu$ with the property that, for any $L>1$, there exists a positive number $\alpha=\alpha_L$ such that  
  \begin{equation}
 \label{existsHaar} \qquad \qquad \qquad \mu  \stackrel{\alpha}{\approx} f_*  \mu, \qquad\qquad \text{ for all } L \text{-biLipschitz maps } f:X\to X .
 \end{equation}
[Uniqueness]  
If moreover $(X,d)$ is $L$-\bilip homogeneous, 
then, whenever another Radon measure $\nu$ also satisfies (\ref{existsHaar}), we have that 
   $\mu\stackrel{\beta}{\approx} \nu$, for some $\beta>1$.
  \end{proposition}
   Measures satisfying  (\ref{existsHaar})  are  called Haar-like. In section 4, we discuss some connections between the existence of a Poincar\'e inequality and upper bounds on the Hausdorff dimension, cf.~Proposition \ref{upperbounddimension}.
 We also show that every Haar-like measure satisfies a lower
and an upper polynomial bound for the measure of balls in terms of the
radius of the ball, cf.~Corollary \ref{growthofballs}.

Before summarizing the strategy for proving Theorem \ref{teorema1}, let us recall some terminology; a standard reference is \cite{Surlesgroupeshyperboliques}.
A geodesic triangle is said to be $\delta$-\textit{thin} if each edge is in the $\delta$-neighborhood of the other two edges. If every geodesic triangle is $\delta$-thin, the space is said to be $ \delta$-\textit{hyperbolic}. A triangle that is not $\delta$-thin is said $\delta$-\textit{fat}.

Here is the intuition behind the proof of Theorem \ref{teorema1}: using charts and a preliminary  argument, cf.~Lemma \ref{unifhomog}, we may suppose that our space is a neighborhood $U$ of the origin $O$ in the plane  $\R^2$ that is uniformly biLipschitz  homogeneous, say with constant $L$.
Then we consider two complementary situations, one is going to imply the theorem, the other will result in  a contradiction.
\begin{description}
\item[ Either] there exists some $\rho$ such that, for any  $r$ smaller than $\rho$,
            there exists an $r/M$-fat triangle in $B(O,r)$;
            $M$ will be a fixed number depending only on the \bilip constant $L$.

            In this case, cf.~Corollary \ref{surrounding triangle}, there exists an $r(10M)^{-1}$-ball
            surrounded by the triangle.
            The basic idea of the argument is to consider the surrounding function $\sur(p,r)$ which is the minimum length of loops that surround the metric ball $B(p,r)$, remind that $B(p,r)$ is now a subset of $\R^2$.
            Therefore,
            the surrounding function for the above ball is less than the length of the triangle's edges, which is less than  $6 r$.
            Using ``quasi invariance'' of the function, cf.~Lemma \ref{qi}, in Corollary  \ref{surbound} we get the existence of some constant $k$ such that, for some $\rho'>0$,
            $$\qquad\sur(p,r) < k r,\quad \forall p\in U,\forall r<\rho'.$$
            From this last bound, we deduce the local doubling and  locally linearly contractible properties, cf.~Proposition \ref{boundedcase} and Proposition \ref{LLC} respectively.
\item[Or]  for any natural number $n$, there exists $r_n < 1/n$ such that
            any triangle in $B(O,r_n)$ is $r_n/M$-thin.

In other words,
             $B(O,r_n)$ is $r_n/M$-hyperbolic.
            BiLipschitz homogeneity implies that any $r_n/L$ ball is $Lr_n/M$-hyperbolic.
            However, such a local hyperbolicity implies, via a corollary of Gromov's coarse version of the Cartan-Hadamard Theorem, cf.~Corollary \ref{gromov}, 
             that the space is globally hyperbolic,
            if we chose $M$ carefully. Set $M=CL^2$ (the constant $C$ is the universal constant in the theorem of Gromov): then
            Gromov Theorem holds  and so
            our initial neighborhood $U$ is $C''r_n$-hyperbolic for any $n\in\N$ ($C''$ is depending only on $L$).
            Since $r_n$ goes to $0$, we have that $U$ is $0$-hyperbolic.
            Every $0$-hyperbolic space  is a tree or an $\R$-tree, cf.~\cite[page 31]{Surlesgroupeshyperboliques}.
            This is a topological contradiction since $U$ is an open set of the plane. This second situation could not in fact occur.
\end{description}

The idea of the construction of Haar-like measures is as follows.  For each $r > 0$,
we consider a maximal $r$-separated net $N_r$.  Let $\mu_r$ be
a sum of Dirac masses at the elements of $N_r$, and re-scale
the result so that the mass of some unit ball $B(x_0,1)$ is $1$.   
Then we claim that the measures $\mu_r$ sub-converge weakly
to the good measure on $X$.  
Now, the existence of the measure is assured by the doubling property and does not require \bilip homogeneity, cf.~Proposition \ref{measure Prop 2}.
The equivalence class of such measures is unique when the space is \bilip homogeneous, cf.~Proposition \ref{measure Prop 1}.

\setcounter{tocdepth}{2}

\tableofcontents
Many thanks go to Bruce Kleiner for inspirational advice and encouragement during the investigation of this problem.
 
\section{Preliminaries}
Throughout all paper, $(X,d)$ will be a {\em locally  biLipschitz  homogeneous} metric space, i.e.,   with the property
that,  
 for every two points
$x_1,x_2 \in X$, there is a pointed biLipschitz homeomorphism $f: (U_1,x_1)\ra (U_2,x_2)$, where $U_i$ is a neighborhood of $x_i$, for $i=1,2$.

\subsection{Uniform \bilip homogeneity}


Given a family $\F$ of homeomorphisms of $X$, we say that $\F$ {\em is transitive on a subset} $U\subset X$ if, for each pair of  points $p,q\in U$, there exists a map $f\in \F$ such that $f(p)=q$.

We will now prove that locally  biLipschitz  homogeneity implies that some 
 family of uniformly \bilip maps, defined on  some neighborhood $U$ of some point,  is transitive  on   $U$. 
 Such argument is based on   Baire Category Theorem and has been used several times in the theory of homogeneous compacta, e.g., in  \cite[Theorem 3.1]{MacManus} or   \cite[Theorem 6.1]{Hohti}.

\begin{lemma} \label{unifhomog} Let $(X, d)$ be any locally compact metric space. 
Suppose $(X,d)$ is  locally   biLipschitz  homogeneous.
Then, for any point of $X$, there exist a compact neighborhood $U$ of the point and   a constant $L$ with the property  that  the family $L$-BiLip$(U;X)$, i.e., the maps  defined on $U$ with values on $X$ that are $L$-biLipschitz,  is  transitive on $U$.
\end{lemma}

\proof Fix a base point $O\in X$ that we will call origin, and let 
$W$ be a compact neighborhood of the origin. Consider the sets
$$S_{n,m}:=\left\lbrace p\in W \;|\;f(O)=p, \;{\rm \;for \;some } \;f:\overline{B\left(O,\frac{1}{m}\right)}\to X,\; n{\rm-\bilip}\right\rbrace .$$
By transitivity, we have $W=\bigcup_{m,n\in\N}S_{n,m}.$
We claim that each $S_{n,m}$ is closed. Take a sequence $p_j\in S_{n,m}$ converging to $p\in W$. Each $p_j$ gives a function $f_j:  \overline{ B(O,\frac{1}{m})}\ra X$. The $f_j$'s are $n$-\bilip, and $f_j(O)=p_j$ converges.
The Ascoli-Arzel\`a argument implies that $f_j$ converges to some $f$ uniformly on the closed ball $\overline{B(O,\frac{1}{m})}$, and the limit function is $n$-biLipschitz.
Therefore, $f(O)=p$ for an $n$-\bilip map $f$ on $\overline{B(O,\frac{1}{m})}$. Thus $p\in S_{n,m}$ and so $S_{n,m}$ is closed.

 Baire Category Theorem implies that there exists an $S_{N,M}$ that has non-empty interior. Therefore  $S_{N,M}$ is a compact  neighborhood of some point $q$. 
Let $f_q:B(O,\frac{1}{M})\ra X$ be an $N$-\bilip map such that $f_q(O)=q$.

We claim that $U:=f_q^{-1}(S_{N,M})\cap B(O,\frac{1}{2MN^4})$ is a neighborhood satisfying the conclusion of the lemma with $L:=N^4$. Indeed, for any two points $p_1, p_2\in f_q^{-1}(S_{N,M})$,
for $i=1,2$, $f_q(p_i)\in S_{N,M}$; so there exists 
an $N$-\bilip map
 $f_i:B(O,\frac{1}{M})\ra X$,  such that $f_i(0)=f_q(p_i)$. Thus we have 
$$p_2=\left( f_q^{-1}\circ f_2\circ f^{-1}_1\circ f_q\right) (p_1)$$
and $f_q^{-1}\circ f_2\circ f^{-1}_1\circ f_q$ is $L$-\bilip. If moreover $p_1\in  B(O,\frac{1}{2MN^4})$, the function  is defined in all $B(O,\frac{1}{2MN^4})$.
$$ \xymatrix{B(O,\frac{1}{2MN^4})\subset B(p_1,\frac{1}{MN^4}) \ar[r]^{\quad\quad\;f_q}&B(f_q(p_1),\frac{1}{MN^3}) \ar[rd]^{f_1^{-1}}&   \\
&	&B(O,\frac{1}{MN^2})\ar[dl]^{f_2}\\
B(p_2,\frac{1}{M})&B(f_q(p_2),\frac{1}{MN})\ar[l]^{f_q^{-1}}&}$$
\qed

\subsection{Existence of fat triangles and Gromov's coarse version of  Cartan-Hadamard Theorem}
From now on, $(X,d)$ will be a {\it biLipschitz homogeneous geodesic surface} as in the assumptions of  Theorem \ref{teorema1}, i.e.,  a geodesic metric space that is topologically equivalent to a surface and is locally  biLipschitz  homogeneous.

The plan for proving Theorem \ref{teorema1} has been  sketched in the introduction. We now proceed in showing the details. 
By Lemma \ref{unifhomog}, and since $X$ is a surface, we have that $X$ is locally isometric to a compact neighborhood $U$ of the origin $O$ in the plane $\R^2$ equipped with a geodesic distance such that, 
for some $L>1$, the action of the $L$-\bilip maps on $U$
is transitive.

We proceed now with the proof that in $U$ there are triangles that are sufficiently fat, i.e., not thin in the sense of Gromov.

\begin{proposition}\label{fat:triangles}
Let $U$ be a  neighborhood of a point $O$ in a geodesic surface $(X,d)$.
Assume that $L$-BiLip$(U;X)$  is  transitive on $U$.
 Then there exist positive constants $M$ and $\rho$ such that  for any  $r<\rho$
            there exists an $r/M$-fat triangle in $B(O,r)$.
 \end{proposition}

The argument for the above proposition will be by contradiction and will be based on Gromov's generalization of Cartan-Hadamard Theorem.
This result states a  
local-to-global phenomenon: if small balls are $\delta$-hyperbolic then the space is $\delta'$-hyperbolic.
The general version of the theorem is the following. 
\begin{theorem}[{Cf. \cite{Gromov-hypgps}, \cite[Theorem 8.1.2]{Bowditch-notes}}]\label{CHT}
There are  constants $d_0$, $C_1$, $C_2$, and $C_3$ with the following property.
Let $X$ be a metric space of bounded geometry. Assume that for some $\de$,
and $d\geq \max(C_1\de,d_0)$, every ball of radius $C_2d$ in $X$
is $\de$-hyperbolic, and the $d$-Rips complex\footnote{The
$d$-{\em Rips complex} $\rips_d(Z)$ of a metric space $Z$ is defined to be the
simplicial complex whose vertex set is $Z$, where distinct points
$x_0,...,x_n\in Z$ span an $n$-simplex in $\rips_d(Z)$ if and only if
$d(x_i,x_j)\le d$ for all $0\le i,j\le n$. } 
 $\rips_d(X)$ is
$1$-connected.  Then $X$ is $C_3d$-hyperbolic.
\end{theorem}
For an exposition of the above theorem, together with the cited definitions,   refer  to the appendix by M. Kapovich and B. Kleiner in \cite{OOS}. 
What we need is the following immediate consequence of Theorem \ref{CHT}. 
\begin{corollary}\label{gromov}
There are  constants $C$ and $C'$ with the following property.
If $X$ is a  simply-connected geodesic metric space such that, for some $R>0$, every 
ball of radius  $ CR$ is $R$-hyperbolic, then 
the space $X$ is $C'R$-hyperbolic.
\end{corollary}


\proof[Proof of Proposition \ref{fat:triangles}]
The idea is to locally use  Corollary \ref{gromov}.
We may assume that $U$ is a simply connected planar domain which we will consider   as a subset of $\R^2$. 
One can show, cf.~ \cite{LeDonne4}, that  there is a subset $A\Subset U$ that has non empty interior and is geodetically closed, i.e., it is a geodesic space.
Clearly, we may assume that $A$ is simply connected, otherwise we add to it all the components  of $U\setminus A$ non containing $\partial U$, and the set would still be geodetically closed. Since $A$ is  a  simply-connected geodesic metric space we can apply Corollary \ref{gromov}.



Let $C$ and $C'$ be the constant in Corollary \ref{gromov}. Set $M=CL^2$. If the conclusion of the proposition  were not true, then, for any natural number $n$, there exists $r_n < 1/n$ such that
            any triangle in $B(O,r_n)$ is $r_n/M$-thin.
In other words,
             $B(O,r_n)$ is $r_n/M$-hyperbolic. We now use $L$-\bilip homogeneity to conclude that any $r_n/L$-ball of $U$ (and so of $A$) is $Lr_n/M$-hyperbolic. 
   Using the definition of $M$, we have that every 
$C\left(\dfrac{r_n}{CL}\right)$-ball is $\dfrac{r_n}{CL}$-hyperbolic.
             Therefore, by Corollary \ref{gromov}, the whole set $A$ is $C'\left(\dfrac{r_n}{CL}\right)$-hyperbolic, for any $n\in\N$.
            Since $r_n$ goes to $0$, this says that $A$ is $0$-hyperbolic.
            Every $0$-hyperbolic space  is a tree or an $\R$-tree.
            This is a contradiction since $A$ is a  set of the plane with non-empty interior. \qed

\subsection{Existence of surrounded balls}
As anticipated in the introduction, as soon as we have a fat triangle, we are interested in looking at a ball inside the triangle of radius proportional to the fatness of the triangle. The purpose is to have a ball surrounded by the triangle.

\begin{definition} \label{surfunction}
A loop $\gamma\subset X$ {\em surrounds} a subset $\Si\subset X$
if $\gamma\cap \Si=\emptyset$ and $\gamma$ separates
$\Si$ from infinity, i.e., any proper path $\R_+\to X$ starting at $\Si$
intersects $\gamma$. In other words, each path in $X$ starting at a point in $\Si$ that escapes every compact set must intersect $\gamma$.
\end{definition}

The ideas in this subsection were partially inspired by   \cite{Papasoglu}.
\begin{proposition}
Let $p$ and $q$ be two points in a geodesic planar and simply connected domain. Let $\gamma$ be a geodesic from $p$ to $q$ and let $\eta$ be another curve from $p$ to
$q$. Suppose $\gamma$ is not contained in the $R$-neighborhood of $\eta$. Then there exists an $R/2$-ball surrounded by
$\gamma\cup\eta$. 
\end{proposition}

\begin{figure}
\centering
\subfigure[The set $V$ is the collection of points at distance less than $r$ from $\eta$ that are `inside' the closed curve $\eta\cup\gamma$.] 
{
    \label{fig:curvnbhd}
    {\epsfysize=1.85truein{\epsfbox{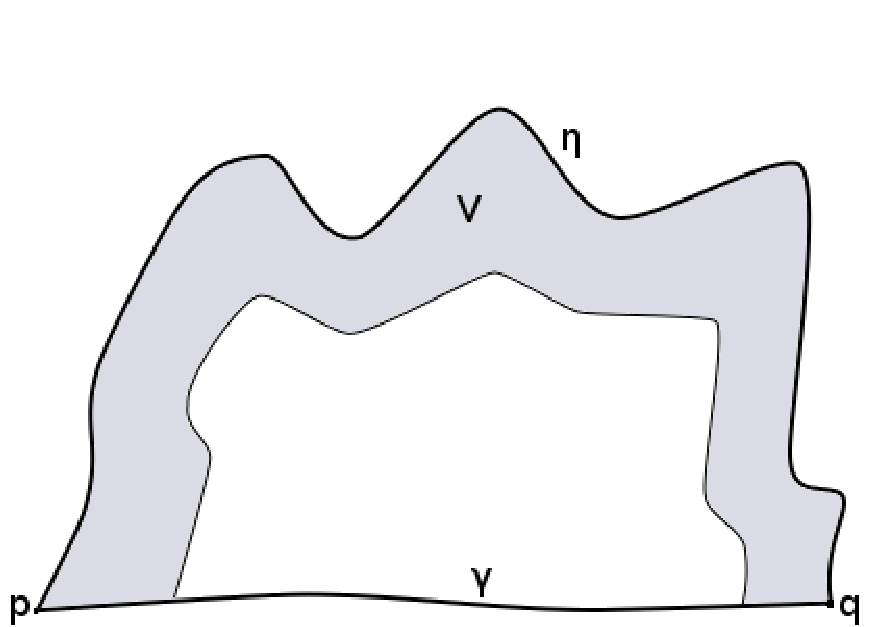}  }}
}
 \hspace{1cm}
\subfigure[An $r$-net in $\sigma$ can be projected on $\gamma$.] 
{
    \label{fig:curvproj}
    {\epsfysize=1.7truein{\epsfbox{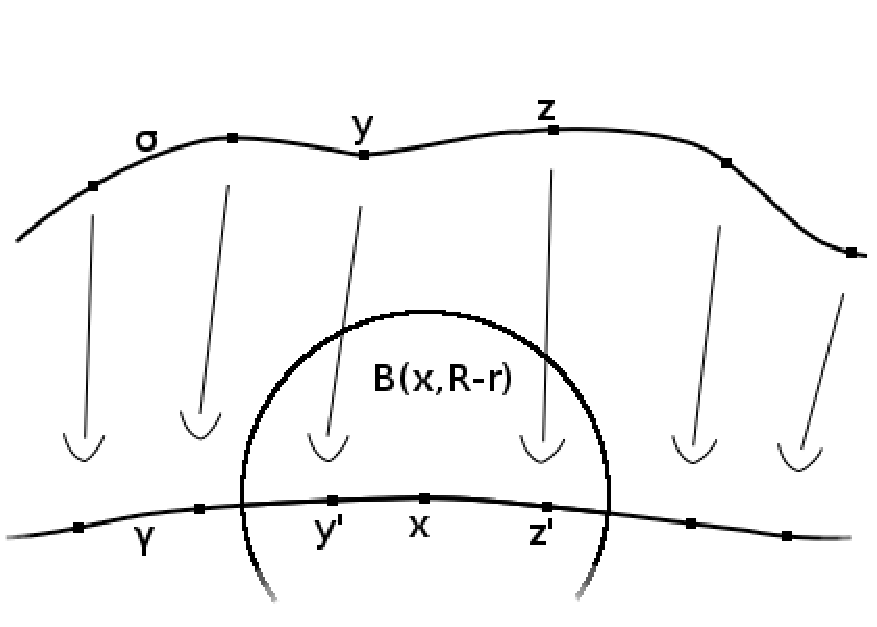}  }}
}
\caption{The proof of the existence of  surrounded balls.}
\end{figure}

\proof Let $0<r<R$. 
Call $U$ the `inside'  $r$-neighborhood of $\gamma$ and $V$ the `inside'  $r$-neighborhood of $\eta$. The word `inside' means that we consider the intersections of the $r$-neighborhoods of the curves with the union of the bounded components of $\R^2\setminus( \eta\cup\gamma)$,
see Figure \ref{fig:curvnbhd}.
 We claim that the complement of $U\cup V$ has a bounded component.
As a consequence,  
 the $r$-ball centered at any point in that component would be surrounded by $\gamma\cup\eta$, and the proof would be concluded.
 If such a claim were not true, then (by Jordan separation) the union $U\cup V$ would be simply connected. Since both $U$ and $V$ are connected,   Mayer-Vietoris Theorem tells us that the intersection $U\cap V$ is connected as well. (Note that $p$ and $q$ are in $U\cap V$). 
Let $\sigma$ be a curve from $p$ to $q$ inside $U\cap V$.
From the hypothesis we know that there exists a ball of radius $R$ and center at some point $x\in\gamma$ that do not intersect $\eta$.

 We claim that $\sigma$ cannot avoid  the ball of center $x$ and radius $R-r$. Assume otherwise. Take an $r$-net along the curve $\sigma$. 
 To each point in the net we can associate a point on $\gamma$, different from $x$, at distance less than $r$, see Figure \ref{fig:curvproj}. The association can be done since $\sigma$ is in the $r$-neighborhood of $\gamma$. This association has to `change sides' of $x$ at some point, 
 in the sense that 
  there are two consecutive points $y$ and $z$ of the net that have associated points $y'$ and $z'$ in disjoint component of $\gamma\setminus \{x\}$, see Figure \ref{fig:curvproj}. 
Now, since both $y$ and $z$ are outside the $(R-r)$-ball,
$$d(x,y')\geq d(x,y)-d(y,y')\geq(R-r)-r=R,$$
and similarly $d(x,z')\geq R$.
This tells us that on one hand, since $y',x,z'$ are on a geodesic in this order,   we have $$d(y',z')=d(y',x)+d(x,z')\geq 2R.$$
On the other hand,
$$d(y',z')\leq d(y',y)+d(y,z)+d(z,z')\leq 3r.$$
But if we choose $r=R/2$, then we get $2R\leq 3r$, which is false.

Thus $\sigma$ intersects the ball of radius $R-r$ center at $x$. 
However, each point in
$\sigma$ is not farther than $r$ from $\eta$. This would imply that $x$ is at distance strictly less than $R$ from $\eta$. This is a contradiction. \qed

\begin{corollary}\label{surrounding triangle}
In a geodesic, planar, and simply connected domain, each geodesic triangle that is not $R$-thin surrounds an $R/2$ ball.\end{corollary}

\proof Let  $\gamma$ be the geodesic edge that is not in the
$R$-neighborhood of the other two edges and let $\eta$ be the concatenation of the other two
edges. Now use the previous proposition. \qed

\subsection{Existence of cutting-through \bilip segments}
\label{segments}
By Lemma \ref{unifhomog}, each point in the space $X$ of Theorem \ref{teorema1} has a neighborhood that is uniformly \bilip homogeneous.
\begin{lemma}\label{cutting-through}
Let $U$ be a  neighborhood of a point $O$ in a geodesic surface $(X,d)$.
Assume that $L$-BiLip$(U;X)$  is  transitive on $U$.
 Then there is a smaller neighborhood $V\subset U$ of $O$ such that, for any $p\in V$,
there is a $L^2$-\bilip image of an interval into $X$ passing through $p$ and starting and ending outside $V$. 
\end{lemma}
\proof 
Take a geodesic $\eta$ in $U$ starting at $O$ and ending at some point $\tilde q$.
Take $V\subset U$ to be a  neighborhood of $O$ such that 
 $$\diam(V)<\dfrac{d(O,\tilde q)}{2L^2}.$$
Let $q\in\eta$ be the midpoint, i.e., 
$$d(O,q)=\dfrac{1}{2} d(O,\tilde q).$$
Take $f$ an $L$-\bilip map such that $f(q)=O$.
For any   point $p\in V$, let $f_p$ an $L$-\bilip map such that $f_p(0)=p$.
Then we claim that $f_p\circ f\circ \eta$ is an $L^2$-\bilip curve passing through $p$ whose end points   are outside $V$.
Indeed,
since $q\in\eta$,
$$p=f_p(O)=f_p(f(q))\in f_p\circ f\circ \eta;$$
  the end point $f_p(f(O))$ lies outside $V$ since
\begin{eqnarray*}
d(f_p(f(O)),p)&=&d(f_p(f(O)), f_p(f(q)))\\
&\geq& \dfrac{1}{L^2}d(O,q)\\
&=&\dfrac{1}{2L^2}d(O,\tilde q)\\
&>&\diam(V)
\end{eqnarray*}
and analogously $d(f_p(f(\tilde q)),p)>\diam(V)$.
\qed

\section{The surrounding function}
We now consider  surrounding functions in 
a  biLipschitz homogeneous geodesic surface  $(X,d)$. 
Studying  linear bounds of surrounding functions is useful for the proof of the doubling property. 
We defined the notion of a loop surrounding a set in Definition \ref {surfunction}. 
If  $\gamma$ is a loop in $X$, we let $ |\ga|$ denote the length of $\gamma$ with respect to the metric $d$.

\begin{definition}[Surrounding function]
\label{surround}
Given $p\in X$, $r\in \R_+$,  let $\sur(p,r)$
be the infimum of lengths of loops
$\gamma \subset X$ that surround the metric ball $B(p,r)\subset X$.
\end{definition}
We actually need a local substitute to control the diameter of the surrounding loops.
\begin{definition}
Given $p\in X$, $r_1,r_2\in \R_+$, with $r_2>r_1$, let $\sur_{r_2}(p,r_1)$
be the infimum of lengths of loops
$\gamma \subset B(p,r_2)$ that surround the metric ball $B(p,r_1)\subset X$.
\end{definition}

Note that, since $X$ is  locally compact, if the set of such loops is non-empty, then there exists a minimum by Ascoli-Arzel\`a Theorem.
We will refer to a loop $\gamma$ that realizes the minimum as a {\em smallest} or {\em shortest} loop that surrounds $B(p,r_1)$.

\begin{lemma}
\label{qi}
The function $\sur_{(\cdot)}(\cdot,\cdot)$ is ``quasi-invariant'': if $f:U\subset X\ra X$ is an $L$-\bilip map with $f(p')=p$ such that $B(p',r_2)\subseteq U$, then
$$\frac{1}{L}\sur_{Lr_2}\left(p,\frac{r_1}{L}\right)\leq\sur_{r_2}(p',r_1).$$ 
\end{lemma}

\proof We may suppose that $\sur_{r_2}(p',r_1)$ is finite.
Choose a smallest loop $\gamma\subset B(p',r_2)$
that surrounds $B(p',r_1)$.  
Then $f(\gamma)$ is a loop that surrounds $B(p,r_1/L)$, is in  
$B(p,Lr_2)$,
and its length is no more than $L|\ga|$. Therefore $\sur_{Lr_2}(p,r_1/L)\leq L|\ga|=L\sur_{r_2}(p',r_1)$.
\qed
 

By Proposition \ref{fat:triangles}, we can now prove the upper bound for the surrounding function.
\begin{corollary}\label{surbound} 
Let $U$ be a  neighborhood of a point $O$ in a geodesic surface $(X,d)$.
Assume that $U$ is homeomorphic to a planar compact domain and that $L$-BiLip$(U;X)$  is  transitive on $U$. 
Then there exist constants $C$ and $\rho'$
such that
            $$\sur(p,r) < C r,$$     for any $p\in U$ and $r<\rho'.$
\end{corollary}
\proof
Let $M$ and $\rho$ be the constants from Proposition \ref{fat:triangles}.
Set $$\rho':=\frac{1}{6ML}\min\{\rho,d(O,U^c)\}\qquad\text{and}\qquad C>12ML^2.$$
For any $r<\rho'$, since $2MLr<\rho$, Proposition \ref{fat:triangles} gives the existence of a $2Lr$-fat triangle in $B(0,2MLr)$.
Corollary \ref{surrounding triangle} says that such a triangle surrounds an $Lr$-ball, which we call $B(\tilde p, Lr)$.

By definition of $\rho'$ and the fact that $\tilde p\in B(0,2MLr)$, we have
\begin{equation}\label{ptildeinU}
B(\tilde p, 4MLr)\subset B(0,6MLr)\subset U .
\end{equation}
Furthermore,
since $B(0,2MLr)\subset B(\tilde p, 4MLr)$ and since the length of a triangle can be bounded by three times the diameter of a ball containing it, we have
$$\sur_{4MLr}(\tilde p, Lr)\le 12MLr.$$
Take any $p\in U$ and take $f:U\to X$ an $L$-\bilip with $f(\tilde p)=p.$
Thus, using Lemma \ref{qi} together with (\ref{ptildeinU}), we finally have 
$$\sur(p,r)\le\sur_{4ML^2r}(p,r)\le L\sur_{4MLr}(\tilde p, Lr)\le (12ML^2)r<Cr.$$
\qed

We now present some technical preliminaries used later for proving local linear connectedness and the doubling property for $U$ in Proposition \ref{LLC} and \ref{boundedcase} respectively.

\begin{lem}\label{various} Suppose to be in the conclusion of Lemma \ref{cutting-through}. Set $R_0:= d(V,U^c)$. We can suppose $R_0>0$. Let $r\in(0, R_0)$ and $p\in V$. Let $\gamma\subset U$  be a loop that surrounds the ball $B(p,r)$.
For the constants $C_0:=2/L^4$, $C_1:=2L^4 $, $C_2:=4L^4$ the following properties are true.

1. We have $\diam(\gamma)\ge C_0 r$.

2.   We have $\sur(p,r)\geq C_0r$.

3. For $r'<C_0r$ and each $p'\in \gamma$, the length of $\gamma\cap B(p', r')$ is at least $r'$.

4. The loop $\gamma$  must lie in $B(p,C_1|\gamma|)$.

5.   
The connected component of $p$ in $X\setminus \gamma$ is contained in $B(p,C_2|\gamma|)$.

\end{lem}

\proof 1. 
Set $K=L^2$. By Lemma \ref{cutting-through} we can consider a $K$-\bilip segment $\sigma$ through $V$ with $\sigma(0)=p$.
Since $\gamma$ surrounds $B(p, r)$ there are two points $p_\pm\in \gamma$
such that $\sigma(T_{\pm})=p_\pm$, with $T_-<0< T_+$. Thus
$$
\diam(\gamma)\ge d(p_-, p_+)\ge \frac{1}{K}(T_+-T_-)  
=\frac{1}{K} [d(T_-,0)+d(0,T_+)]$$
$$
\ge\frac{1}{K^2} d(p_-, p)+ d(p, p_+)
 \ge \frac{1}{K^2} (r+r)= \frac{2}{K^2} r.
$$

2. Let $\gamma$ be a smallest loop surrounding $B(p, r)$. By part (1) 
$$
\sur(p, r)=|\gamma|\ge \diam(\gamma) \ge \frac{2}{K^2} r.
$$

3. According to (1), $\diam(\gamma)\ge C_0 r$.
Hence, for each $r'\leq C_0r$    and $p'\in \gamma$, the metric sphere
$S(p', r')$ has nonempty
intersection with $\gamma$. Thus the length of $\gamma\cap B(p', r')$
is at least $r'$.

4.   
Let $p_\pm$ be the points considered in (1). Then 
$$d(p_\pm, p)\le Kd(T_\pm,0)=K|T_\pm|  \le K |T_+- T_-|\le K^2 d(p_+, p_-) \le K^2|\gamma|.$$
Thus, for any  $z\in\gamma$,
$$
d(z, p)\le d(z, p_+)+d(p_+,p)\le|\gamma|+K^2|\gamma|\le 2 K^2|\gamma|.
$$ 
Therefore $\gamma\subset B(p, C_1|\gamma|)$.

5. Consider a point $q\in X\setminus \gamma$ that lies in the same component of
$X\setminus\gamma$ as $p$. Then either $q\notin {\rm Nbhd}_r(\gamma)$ or $d(q, \gamma)\le r$.
In the first case the loop $\gamma$ surrounds both $B(p,r)$ and $B(q,r)$. Hence, by (4)
$$
\gamma\subset B(p, C_1|\gamma|)\quad\text{ and }\quad \gamma\subset B(q, C_1|\gamma|),
$$
i.e., any point of $\gamma$ is at distance less than $C_1|\gamma|$ from both $p$ and $q$.
By the triangle inequality we conclude that 
$d(p, q)\le 2C_1|\gamma|=4K^2|\gamma|$.
In the second case, if $d(q, \gamma)\le r$, then $d(q, \gamma(t))\le r$, for some $t$.
Thus, using (4) and (2), we have $$d(p, q)\le d(p,\gamma(t))+d(\gamma(t),q)\le C_1|\gamma|+r\le C_1|\gamma|+|\gamma|/C_0=\frac{5}{2}K^2|\gamma|<4K^2|\gamma|.$$ Therefore
$q\in B(p, C_2|\gamma|)$. 
\qed

\subsection{Local linear contractibility and the doubling property}
We are ready to prove local linear contractibility and the local doubling property,  
for \bilip homogeneous geodesic surface. 

\begin{proposition}
\label{LLC} 
Let $(X,d)$ be a biLipschitz homogeneous geodesic surface. 
Then any point of $X$ has a  locally linearly contractible neighborhood.
\end{proposition}
\proof
Fix a point $O\in X$. Let $U$ be the neighborhood given by Lemma \ref{unifhomog},
so we are in the assumption of Corollary \ref{surbound}. Let $V$ be the  
 neighborhood for $O$  given by Lemma \ref{cutting-through}. Therefore the conclusions of Lemma \ref{various} hold. 
We will consider $U$ and $V$ as planar domain of $\R^2$.  On them 
the surrounding function has the linear upper bound, by Corollary \ref{surbound}.

Let $p\in V$ and  $r<\rho'$, so that Corollary \ref{surbound} holds.
Consider the ball $B(p,r)$ and a length minimizing surrounding loop $\gamma$. 
Note that the ball $B(p,r)$ is connected, being the metric geodesic.
Since $\gamma$ is minimizing, then it is a simple loop. Thus, by Jordan Theorem,  the bounded component $E$ of $\R^2\setminus \gamma$ is homeomorphic to a disk, in particular $E$ is homotopic to a point. 
Since the ball $B(p,r)$ is connected, it is contained in $E$. Thus
 $B(p,r)$ is homotopic to a point in $E$. By point (5) in the previous lemma, $E$ is contained in $B(p,C_2|\gamma|)$.
The bound on the surrounding function gives $|\gamma|=\sur(p,r)<Cr$ and so $B(p,C_2|\gamma|)\subset B(p,C_2Cr)$.
In conclusion, $B(p,r)$ is homotopic to a point in $B(p,C_2Cr)$. \qed

\begin{proposition}
\label{boundedcase} 
Let $(X,d)$ be a biLipschitz homogeneous geodesic surface. 
Then any point of $X$ has a   neighborhood that is doubling.
\end{proposition}
\proof  
As in the previous proof, 
fixed a point $O\in X$, let $U$ and $V$ be the neighborhoods given by Lemma \ref{unifhomog} and Lemma \ref{cutting-through} respectively. 
We will consider $U$ and $V$ as planar domain of $\R^2$. 
Thus the conclusions of Corollary \ref{surbound}  and  Lemma \ref{various} hold.
In particular,  we have the upper bound for the surrounding function.
Namely, for $p\in V$ and  $r<\rho'$, if $\gamma$ surrounds $B(p,r)$ and is a minimizer for $\sur(p,r)$, then, by Corollary \ref{surbound}, we have $|\ga|\leq Cr$.
Moreover, Part (5) of Lemma \ref{various} says that   the connected component of $p$ in $X\setminus \gamma$ is contained in $B(p,C_2Cr)$, since $C_2Cr\geq C_2|\ga|$.

Fix $p\in V$. Choose a loop $\gamma_1\subset X$ with length at most $Cr$ that surrounds $B(p,r)$, and set $\L_1=\{\gamma_1\}$.
Let $\n_1$ be an $\frac{r}{2L^2}$-separated $\frac{r}{2}$-net in $\gamma_1$.  Then, by Lemma \ref{various} (3),
the cardinality of $\n_1$ is at most
$$
\frac{|\gamma_1|}{r/(4L^4)}\le \frac{Cr}{r/(4L^4)}=4L^4C=:c.
$$
Let $\L_2$ be a collection of loops (each having length at most $Cr$) surrounding the $r$-balls centered at the points in $\n_1$.
Proceed inductively in this fashion, building up $k$ layers of surrounding loops in $X$.
The union $\mathcal V_k\defeq\n_0\cup\ldots\cup \n_k$ has cardinality at most
$$
c^{k+1}=(4L^4C)^{k+1}.
$$
We claim that the collection of $C_2Cr$-balls centered at the points in $\mathcal V_k$ covers $B(p,\frac{kr}{2})$.  To show such a claim, consider a path $\si$ of length at most $\frac{kr}{2}$ starting at $p$.
Inductively break $\si$ into a concatenation of at most $k$ sub-paths of length at least $\frac{r}{2}$ as follows.  Let $\si_1$ be the initial segment of $\si$ until $\si$ intersects  $\gamma_1$. The path $\si_1$ has length at least $r$ and terminates within distance $\frac{r}{2}$ of a point $p_1\in\n_1$. Let $\si_2$ be the initial segment of $\si\setminus\si_1$ until $\si\setminus\si_1$ intersects  the surrounding loop for  $B(p_1,r)$,
et cetera.  At each step the segment $\si_i$ has length at least $\frac{r}{2}$, and from what was said at the beginning of the proof, each $\si_i$ is contained in $\cup_{q\in  \mathcal V_k}B(q,C_2Cr)$.

Thus
$$
B\left(p,\frac{kr}{2}\right)\subset \bigcup_{i=1}^{c^{k+1}} B(p_i, C_2Cr),
$$
for each  $p\in V$.
Choose $k$ such that $\frac{k}{2}= 2C_2C$ (clearly we may assume $C_2, C\in \N$) and 
define the constant $N:= c^{k+1}$.
Writing $\rho$ in the form $\rho=C_2Cr$,
we have proved that, for any $p\in V$,
$$
B(p, 2\rho)\subset \bigcup_{i=1}^{N} B(p_i, \rho).
$$
In other words  $V$ is doubling.
\qed 


\section{Consequences of the doubling property}
\subsection{Dimension consequences}
Recall that doubling spaces are precisely those  spaces with {\em finite Assouad dimension} (also known as metric covering dimension or uniform metric dimension in the literature). See   \cite{Heinonenbook} for the definition. However, the Assouad dimension of a metric space can be defined equivalently as the infimum of all numbers $D>0$ with the property that every ball of radius $r>0$ has at most $C\eps^{-D}$ disjoint points of mutual distance at least $\eps r$, for some $C\geq 1$ independent of the ball.

Let us recall that a set $N\subset X$ is said to be $\epsilon$-{\em separated} if $d(x,y)\geq\epsilon$ for each distinct $x,y\in N$. 
Also, a set $N\subset X$ is said to be an 
 $\epsilon$-{\em net} if, for each $x\in X$, $d(x,N)\le\epsilon$. 
 Clearly an $\epsilon$-separated set that is maximal with respect to inclusions of sets id an $\epsilon$-net;    such a set is called a maximal $\epsilon$-separated net. 
 
 Thus, a metric space $X$ of Assouad dimension less than $D$ has the property that there exists a constant $C$ such that, for any $p\in X$ and  any $r>0$, 
\begin{equation}\label{Assouad}
N_\delta \text{ is }  \delta \text{-separated }\;\Longrightarrow\;\#((N_\delta \cap B (p, r )) \leq  C \left(\dfrac{\delta }{r} \right )^{-D} .
\end{equation}

Since the Hausdorff dimension of a metric space does not exceed its Assouad dimension, the next corollary is immediate.
\begin{corollary}
A locally biLipschitz homogeneous geodesic surfaces has finite Hausdorff dimension.
\end{corollary}
\proof By Proposition \ref{boundedcase} any point has a neighborhood that is doubling. Thus the Hausdorff dimension of such a neighborhood is finite, say $\alpha$. Now, since the space is biLipschitz homogeneous and \bilip maps preserve Hausdorff dimension, all points have neighborhoods with Hausdorff dimension equal to $\alpha$. Since the Hausdorff dimension depends on local data, the dimension of the space is $\alpha$.
\qed

\subsection{Good measure class: the Haar-like measures}
We will now give   the details regarding the   Haar-like  measures. Throughout this section, let $O$ be a fixed point and let $B_r=B_r(O)$. 
Let $\delta_p$ be the Dirac measure defined by $\delta_p(A)=1$ if $p\in A$ and $\delta_p(A)=0$ if $p\notin A$. 
   
\subsubsection*{Notation $\mu  \stackrel{\alpha}{\approx}   \nu$} For  $\mu$ and $\nu$   Borel measures and a number $\alpha>0$, we say that $\mu  \stackrel{\alpha}{\approx}   \nu$ if    $$\frac{1}{\alpha} \nu(A)\le \mu(A)\le\alpha  \nu(A),$$ for each Borel set $A$. Equivalently, if they are absolute continuous with respect to each other and the Radon-Nikodym derivatives are bounded between $\frac{1}{\alpha} $ and $\alpha $, i.e., there exists a function $h:X\to[\frac{1}{\alpha} ,\alpha ]$ so that $d\nu=hd\mu$. 
 
  For a set  $N\subset X$ such that $ \# (N\cap B_1)<\infty$, define the Radon  measure 
  $$  \mu_N := \frac{1}{\# (N\cap B_1)}\sum_{p\in N}\delta_p , $$
  i.e.,
  $$  \mu_N (A)= \frac{\#(N\cap A)}{\# (N\cap B_1)},  $$
for any Borel set $A$.
 The normalization has the purpose of having  $  \mu_N (B_1)= 1$, for any set $N$.
 
Now, the existence of a good measure is assured by the doubling property, and does not require homogeneity. 
Recall that, if $f:X\to X$ is any Borel function, then any Borel measure $\mu$ on $X$ can be pushed forward as
$$\left(f_*\mu\right)(A):=\left(f_\#\mu\right)(A):=\mu\left(f^{-1}(A)\right),$$
for any Borel set $A$.
\begin{proposition}[Existence]\label{measure Prop 2}  Let $(X,d)$ be any  doubling metric space. 
 Then there exists a  non-zero Radon measure $\mu$,
 with the property that, for any $L>1$, there is a constant $\alpha=\alpha_L $ such that $\mu  \stackrel{\alpha}{\approx}   f_*\mu$, for each $f\in L{\rm-BiLip}(X,d)$. 
 \end{proposition}
\proof  
  For each $\epsilon>0$ choose a  maximal $\epsilon$-separated net $N_\epsilon$ and consider the associated measure 
 $\mu_\epsilon:=\mu_{N_\epsilon}$ defined as above, i.e.,
  $$  \mu_\eps (A):= \frac{\#(N_\eps \cap A)}{\# (N_\eps\cap B_1)},  $$
 for any Borel set $A$.

 By Theorem 1.59 in \cite{Ambook}, since the $\mu_\epsilon$ are (finite) Radon measures and $  \mu_\epsilon(B_1)= 1$, there is a sub-sequence $\mu_{\epsilon_n}$ that is weak$^*$ convergent to a measure  $  \mu $. 
 Recall that, if $cl(B)$ is the closure of a set, then
\begin{eqnarray}\label{limit measure}
 \limsup   \mu_{\epsilon_n}(cl(B)) \le  \mu(cl(B))  .
  \end{eqnarray}
 
 Let us prove that $\mu$ satisfies the conclusion of the proposition. Take any $f\in L\text{-Bilip}(X,d)$. Then note that $f(N_\epsilon) $ is an $\frac{\epsilon}{L}$-separated $L \epsilon$-net.
 
 Fix any ball $B$. Take two  other balls $B'' \subsetneq B' \subsetneq B$ with same center and different radii $r''< r'< r$.
If $\epsilon\le r' - r''$ and
 $B(p,\epsilon)\cap B''\neq\emptyset $, then we have that $p\in  B'$. Thus
 $$ B''\subset \bigcup_ {p\in B'\cap N_\epsilon}B(p,\epsilon),\qquad \forall \epsilon\le r' - r'',$$
 since  $N_\epsilon $ is an $\epsilon$-net.
 Moreover, since $f(N_\epsilon)$ is $\frac{\epsilon}{L}$-separated, from (\ref{Assouad}) we have
 $$\# ( B(p,\epsilon)\cap f(N_\epsilon)) \le C \left( \dfrac{\epsilon/L}{\epsilon}\right)^{-D}= C  L^D.$$
 Then 
$$  \# \left(B''\cap f(N_\epsilon)\right) \le \sum_ {p\in B'\cap N_\epsilon} \# \left(B(p,\epsilon)\cap f(N_\epsilon) \right)\le   C  L^D   \# (B'\cap N_\epsilon) .$$
 So,
\begin{eqnarray*}
f_* \mu_\epsilon (B'') &\le& \frac{\#\left (f^{-1}(B'')\cap N_\epsilon\right)}{\# (B_1\cap N_\epsilon)} \\
&=& \frac{\# \left(B''\cap f(N_\epsilon)\right)}{\# (B_1\cap N_\epsilon)} \\
 &\le& C   L^D  \frac{\# (B'\cap N_\epsilon)}{\# (B_1\cap N_\epsilon)} \\
 &=& C   L^D \mu _\epsilon (B') \\
 &\le& C   L^D \mu _\epsilon (cl(B')).
 \end{eqnarray*}
 Taking the limit for $\epsilon_n\rightarrow 0$, from the last estimate and from (\ref{limit measure}), we have
\begin{eqnarray*}
 f_* \mu (B'') & \le&  \liminf   f_* \mu_{\epsilon_n} (B'') \\
 &\le& \limsup   C   L^D \mu_{\epsilon_n}(cl(B')) \\
 &\le& C   L^D \mu(cl(B')) \\
 & \le & C   L^D \mu (B).
  \end{eqnarray*}
Since   $B''\subset B$ was arbitrary, we get $$ f_* \mu (B) \le  C   L^D \mu (B).$$
In conclusion, $f_* \mu\leq \alpha \mu$, for $\alpha:=C   L^D$, on every (small) ball, so the same inequality holds on every open set and therefore on every Borel set.
Since $f^{-1}\in L\text{-Bilip}(X,d)$, we also get 
 $$\dfrac{1}{\alpha} \mu (A) \le f_*\mu (A),$$
  for each Borel set $A$. So both the required inequalities are proven.
\qed

  The equivalence class of the Haar-like measures is unique when the space is \bilip homogeneous.
\begin{proposition}[Uniqueness]\label{measure Prop 1}  
  Let $(X,d)$ be a doubling metric space with a transitive set $\mathcal F$ of $L$-bilip maps. Suppose that two non-zero Radon measures $\mu_1$ and $\mu_2$
  on $X$ are such that  $\mu_i  \stackrel{\alpha}{\approx} f_*  \mu_i$,  for $i=1,2$ and for each $f\in \mathcal F$. Then 
   $\mu_1\stackrel{\beta}{\approx} \mu_2$, for a constructive $\beta>1$.
  \end{proposition}

Let us prepare for the proof of the uniqueness of the class of good measures with a lemma which will be  useful again   in the proof of polynomial growth of such measures.   The following lemma says that if $\mu$ is a Haar-like measure, then the $\mu$ 
measure of the $\eps$-balls is  approximately  
the inverse of the cardinality of a maximal $\eps$-separated net in the unit ball. 
\begin{lemma}\label{measure of balls}
  Let $(X,d)$ be a doubling metric space with a transitive set $\mathcal F$ of $L$-\bilip maps.
  Suppose that a non-zero Radon measure $\mu$
  on $X$ is such that  $\mu  \stackrel{\alpha}{\approx} f_*  \mu$, for each $f\in \mathcal F$.
Then there are positive constants $\eps_0$, $k$, and $h$ such that, for any $\eps<\eps_0$ and for any maximal $\eps$-separated net $N_\epsilon$, defining
$c_\eps:=\#\left(N_\epsilon\cap B_1\right)$, we have
  \begin{equation}\label{measure great}
  \mu \left( B(p,L  \epsilon)\right)\ge k c_\epsilon ^{-1},
  \end{equation}
and
\begin{equation}\label{measure less}
\mu \left( B\left(p, \frac{\epsilon}{2L}\right)\right)\le  h c_\epsilon ^{-1}.
\end{equation}
\end{lemma}

 \proof
Set $\eps_0=1/2$. Let $\eps<\eps_0$. 
    Fix   $p\in X$. For any $p_j$ in the maximal $\epsilon$-separated net $ N_\epsilon$, choose  $f_j\in \mathcal F$ such that $f_j (p)=p_j$. Thus 
    $B(p_j,\epsilon)\subset f_j \left(B(p,L\epsilon)\right)$.

      To show (\ref{measure great}), consider that, since $N_\epsilon$ is an $\epsilon$-net, the family $\left\{B(p_j,\epsilon)\right\}_{p_j\in N_\epsilon}$ is a cover of $X$. Therefore
    $$B_\frac{1}{2} \subset \bigcup\left\{ B(p_j,\epsilon) :    p_j\in N_\epsilon \cap B_1     \right\},$$
because $\epsilon<\frac{1}{2}$ (we had to reduce to the ball $B_\frac{1}{2}$ since, removing those $\eps$-balls with center outside $B_1$, we might fail to cover $B_1\setminus B_{1-\eps}$).
So   
\begin{eqnarray*}
0<\mu ( B_\frac{1}{2}) &\le& \mu \left(    \bigcup \left\{ B(p_j,\epsilon )\;:\;{p_j\in N_\epsilon \cap B_1 } \right\}\right) \\
&\le&   \sum_{p_j\in N_\epsilon \cap B_1 } \mu \left(B(p_j,\epsilon )\right)\\
&\le &\sum \mu \left(f_j \left(B(p,L \epsilon) \right) \right)\\
&=&\sum \left((f_j ^{-1})_*    \mu \right)\left( B(p,L\epsilon)\right)  \\
&\le&   \sum_{p_j\in N_\epsilon \cap B_1 }  \alpha   \mu \left( B(p,L  \epsilon)\right) \\
&=&\# (N_\epsilon\cap B_1)   \cdot   \alpha  \mu \left( B(p,L  \epsilon)\right)   \\
&=&  c_\epsilon      \alpha    \mu \left( B(p,L\epsilon)\right).  
\end{eqnarray*}
Setting $k=\alpha^{-1}  \mu ( B_\frac{1}{2})$, we obtain (\ref{measure great}).
   
Now we show (\ref{measure less}). Since $N_\epsilon$ is  $\epsilon$-separated and $\eps<1/2$, we have that  $\left\{B(p_j,\frac{\epsilon}{2})\right\}_{p_j\in N_\epsilon \cap  B_1}$ is a    disjoint family of subsets of $B_\frac{3}{2}$.
   Therefore,
  \begin{eqnarray*}
\mu ( B_\frac{3}{2}) &\ge& \mu    \left( \bigcup \left\{ B \left(p_j,\frac{\epsilon}{2} \right)\;:\; {p_j\in N_\epsilon \cap B_1 }\right\} \right) \\
&= &  \sum_{p_j\in N_\epsilon \cap B_1 } \mu \left(B \left(p_j, \frac{\epsilon}{2} \right) \right)\\
&\ge& \sum \mu \left(f_j \left(B \left(p,\frac{\epsilon}{2L}\right)\right) \right)\\
&=&\sum \left((f_j ^{-1})_*     \mu\right) \left( B \left(p,\frac{\epsilon}{2L}\right)\right)  \\
&\ge &          \sum_{p_j\in N_\epsilon \cap B_1 }  \alpha^{-1}     \mu \left( B \left(p,\frac{\epsilon}{2L}\right) \right) \\
&=&\# (N_\epsilon\cap B_1)  \cdot     \alpha^{-1}     \mu \left( B \left(p, \frac{\epsilon}{2L}\right) \right)  \\
& =&  c_\epsilon      \alpha^{-1}     \mu \left ( B \left(p,  \frac{\epsilon}{2L}\right) \right). 
     \end{eqnarray*}
Setting $h=\alpha  \mu ( B_\frac{3}{2}) $, we obtain (\ref{measure less}).\qed

  \proof[Proof of Proposition \ref{measure Prop 1}]   Let $s=h/k  $. Then (\ref{measure great}) and  (\ref{measure less}) imply that, for each $\epsilon<\frac{1}{2}$, we have
\begin{equation}\label{measure 1 less measure 2}
\mu_1 \left( B \left(p, \frac{\epsilon}{2L} \right)\right)\le   s  \mu_2 \left( B \left(p,L  \epsilon \right)\right).
\end{equation}   
   
   Now we plan to estimate the measure $\mu_2 \left( B \left(p,L  \epsilon \right)\right)$ with a constant times
   $\mu_2 \left( B \left(p, \dfrac{\epsilon}{2L} \right)\right)$
    using the fact that $(X,d)$ is doubling.  
 Indeed, there is a number $m$, not depending on $\eps$, so that  $ m$ balls of 
radius $\eps /L$ cover $B(p, L \eps )$. Let $q_1 , q_2,\ldots,q_m\in X$ be such that 
$$B(p, L \eps ) \subset \bigcup^m_{i=1} B(q_i , \eps /L).$$
For each $i=1,\ldots,m$, choose $g_i \in  \F$ with $ g_i (q_i ) = p$. Then
\begin{eqnarray*} 
\mu_2 (B(p, L\eps )) &\leq&  \sum^m_{i=1}   \mu_2 \left(B \left(q_i , \eps /L \right)\right)\\
 &\leq&  \sum^m_{i=1}   \alpha \mu_2 \left(f \left(B \left(q_i , \eps /L\right) \right)\right) \\
 &\leq&  \sum^m_{i=1}   \alpha \mu_2 \left(B \left(p, \eps /2L \right)\right)\\
& =& m\alpha \mu_2 \left(B\left(p, \eps /2L \right)\right).
 \end{eqnarray*}
Hence, from (\ref{measure 1 less measure 2}), we have that there exists $\gamma > 0$, such that, for all $ \eps  > 0$, $$\mu_1 \left(B\left(p,\dfrac{\eps}{2L}   \right)\right) \leq  \gamma  \;\mu_2 \left(B\left(p, \dfrac{\eps}{2L}  \right)\right).$$ 
In conclusion, $\mu_1$  is smaller than $\gamma  \mu_2$ on every small ball, so the same is true on every open 
set and thus on every Borel set. 
The symmetric hypothesis on $\mu_1$ and $ \mu_2$ gives us the other inequality too. \qed

\begin{lemma}\label{HausBound}
Let $(X,d)$ be a metric space where a ball $B_{ {1}/{2}}$  has Hausdorff
dimension $\alpha$. 
Then, for any $t>0$ and $c>0$, there exists an $\eps_0>0$ such that any $\eps$-net $N_\eps$ with $\eps<\eps_0$ has the property that
$$\#\left(N_\epsilon\cap B_1\right)\geq \dfrac{c}{\eps^{\alpha-t}} .$$
\end{lemma}
\proof
Since the Hausdorff dimension is $\alpha$, all the Hausdorff measures of dimension less than $\alpha$ are infinite:
\begin{equation}\label{alpha-s Hausdorff}\H^{\alpha-s}(B_{1/2})=\infty, \qquad\forall s>0.\end{equation}

Let us  assume that the conclusion of the lemma is not true, i.e., there exist $t, c>0$ so that, for all $\eps_0>0$, there is an  $\eps$-net $N_\eps$, with $\eps<\eps_0$ with
\begin{equation}\label{epsseq}
\#\left(N_\epsilon\cap B_1\right)\leq \dfrac{c}{\eps^{\alpha-t}} .
\end{equation}
Since  $N_\eps$ is an $\eps$-net, the collection of balls
$$\left\{ B(p,\eps) \;:\;p\in N_\epsilon\cap B_1\right\},$$
for $\eps<1/2$, is a covering of $B_{1/2}$ by sets of diameter less than $2\eps$. We can estimate the Hausdorff measure
\begin{eqnarray*}
\H^{\alpha-s}_{2\eps}(B_{1/2})&:=&\inf\left\{\sum(\diam V_i)^{\alpha-s}\;:\;  \diam V_i\leq2\eps,B_{1/2}\subset\cup V_i\right\}\\
&\leq& \sum_{ p\in N_\epsilon\cap B_1}\left(\diam B(p,\eps)\right)^{\alpha-s}\\
&\leq& \sum_{ p\in N_\epsilon\cap B_1} (2\eps)^{\alpha-s}\\
&\leq& \#\left(N_\epsilon\cap B_1\right)\cdot (2\eps)^{\alpha-s}\\
&\leq&\dfrac{c}{\eps^{\alpha-t}} \; (2\eps)^{\alpha-s}\\
& =& 2^{\alpha-s}c \eps^{t-s}.
\end{eqnarray*}
Taking $s\in(0,t)$, we have that $2^{\alpha-s}c \eps^{t-s}\to 0$, as $\eps\to0$. 
Thus,
for the infinitesimal sequence of $\eps$'s where \eqref{epsseq} holds, we have that $\H^{\alpha-s}_{2\eps}(B_{1/2})$ goes to zero as well. Therefore
$$\H^{\alpha-s}(B_{1/2}):=\lim_{\delta\to0} \H^{\alpha-s}_{\delta}(B_{1/2})=0,$$
contradicting (\ref{alpha-s Hausdorff}).\qed

Let us remark that since $(X,d)$ is doubling, the cardinality of $N_\epsilon\cap B_1$ is finite. In fact, using (\ref{Assouad}), such a cardinality is bounded by 
    $C\eps^{-D},$ for some constants $C>0$ and any $D$ greater than the Assouad dimension.  
    Using Lemma \ref{HausBound} and Lemma \ref{measure of balls} we conclude the following. 
   
   \begin{corollary}\label{growthofballs}
Let $(X,d)$ be a doubling $L$-\bilip homogeneous metric space. Let $\mu$ be a Haar-like measure. Then, for any  $t>0$, there exists $r_0>0$ and $K>1$ such that, for all $p\in X$ and any $r<r_0$,
   $$\frac{1}{K}\;r^{\dim_{\rm A}(X,d)+t}<\mu\left(B(p,r)\right)<K\;r^{\dim_{\rm H}(X,d)-t}.$$
\end{corollary}

Recall that $\dim_{top}\leq\dim_{H}\leq\dim_{A},$
so for $r<1$, we have
$r^{\dim_{A} }\leq r^{\dim_{H}}\leq r^{\dim_{top}}.$

\begin{corollary}
Let $\gamma$ be a rectifiable curve. For any Haar-like measure $\mu$, we have $\mu(\gamma)=0$.
\end{corollary}

Since any doubling measure is non-atomic and strictly positive on non-empty open sets, we are allowed to use the following theorem by Oxtoby and Ulam.
\begin{theorem}[{\cite[Theorem 2]{Oxtoby-Ulam}}]
Let $\mu$  be a Radon measure on the square $Q=[0,1]^n$, $n>2$, with the properties that
\begin{itemize}
\item[(i)] $\mu$ is zero on points,
\item[(ii)] $\mu$ is strictly positive on non-void open sets,
\item[(iii)]  $\mu(Q)=1$,
\item[(iv)]  $\mu(\partial Q)=0$.
\end{itemize}
Then there exists an homeomorphism $h: Q\to Q$ such that
$\mu=h_* \L$.
\end{theorem}
As an immediate consequence we have the following:
\begin{corollary}
Any doubling measure on the plane is locally a multiple of the Lebesgue measure up to a continuous change of variables.
\end{corollary}

 \subsection{Upper bounds for the Hausdorff dimension}
 It is an open question whether a \bilip homogeneous geodesic surface satisfies a Poincar\'e inequality. However, we now show that the existence of a Poincar\'e inequality implies a bound  on the Hausdorff dimension.
 
 Let $1\leq p<\infty$. We say that a measure metric space $(X,d,\mu)$ admits a weak $(1,p)$-Poincar\'e inequality if there are constants
$\lambda\geq 1$ and $C\leq1$ so that
$$\dashint_B|u-u_B|\;d\mu  \leq C(\diam B)\left( \dashint_{\lambda B}\rho^p \;d\mu\right)^{1/p} ,$$
for all balls $B\subset X$, all bounded continuous functions $u$ on $B$, and all upper gradients $\rho$ of $u$. 
Recall that $\rho$ is an upper gradient for $u$ if
$$|u(x)-u(y)|\leq \int_{\gamma_{xy}}\rho \;ds,$$
for each rectifiable curve $\gamma_{xy}$ joining $x$ and $y$ in $X$.

\begin{proposition}\label{upperbounddimension} Let $(X,d)$ be a  \bilip homogeneous geodesic surface.
If a weak $(1,p)$-Poincar\'e inequality holds for a Haar-like measure $\mu$, then  
$$\dim_{H} (X,d)\leq 1+p.$$
\end{proposition}

\proof
 We may assume that $X$ is a planar domain. Fix any geodesic $\sigma$ in $X$. Consider a simply connected set $B\subset X$ that is divided into two parts by $\sigma$, i.e., $B\setminus\sigma=A_0\sqcup A_1$ with $A_0$ and $A_1$ simply connected. Define the following functions: 
$$\delta(p):=\left\{\begin{array}{rcc} d(p,\sigma)& \text{for} &p\in A_1 \\ -d(p,\sigma)& \text{ for } &p\in A_0 \end{array}\right. \text{ and }
u_\eps(p):=\left\{\begin{array}{ccc} 
\dfrac{\eps-\delta(p)}{2\eps}& \text{for}& -\eps\leq\delta(p)\leq\eps \\
0 &\text{for} &\delta(p)\leq-\eps \\
1 &\text{for} &\delta(p)\geq\eps  \end{array}\right. .
$$
The function $u_\eps$ is $0$ on those points of $A_0$ at distance more than $\eps$ from $\sigma$. In the $\eps$-neighborhood of $\sigma$ it increases linearly in the distance from $\sigma$ to the value $1$ at those   points of $A_1$ at distance more than $\eps$ from $\sigma$. Therefore the function $\rho_\eps$ defined to be $\dfrac{1}{2\eps}$ on the $\eps$-neighborhood of $\sigma$ and $0$ elsewhere is an upper-gradient for $u_\eps$.

Since $u_\eps\to\chi_{A_1}$ as $\eps\to0$, an easy computation gives that 
$$\dashint_B|u_\eps-(u_\eps)_B|\;d\mu \to \dfrac{2\mu(A_0)\mu(A_1)}{(\mu(B))^2}\neq 0.$$
So the limit is non-zero.

Let us now see how  the Poincar\'e inequality estimates the previous limit.
 Cover the $\eps$-neighborhood of $\sigma$ with $\dfrac{{\rm length}(\sigma)}{\eps} $ balls of radius $2\eps$. Thus, if $\alpha$ is any number smaller than the Hausdorff dimension, using Corollary \ref{growthofballs}, we get
\begin{eqnarray*}
\left( \dashint_{\lambda B}\rho^p \;d\mu\right)^{1/p} &\leq & \left( \sum_j(\mu(B(p_j,2\eps)))  (\dfrac{1}{2\eps})^p \right)^{1/p}\\
&\leq & \left(\dfrac{{\rm length}(\sigma)}{\eps} \dfrac{K(2\eps)^\alpha}{(2\eps)^p}\right)^{1/p}=K' (\eps^{\alpha-1-p})^{1/p}.
\end{eqnarray*}
If it would be possible to have $\alpha>1+p$, then this last term would go to zero, as $\eps$ goes to zero, and it would give a contradiction. So $\alpha$ and hence $\dim_{H} (X,d)$ must be smaller than $1+p$. \qed

An immediate consequence of the above proposition is that  the existence of a $(1,1)$-Poincar\'e inequality implies that the Hausdorff dimension is $2$.
\subsection{Lower bound for the Hausdorff $2$-measure}


Another consequence of the lower bound  on the surrounding function is a  density bound on
the $2$-dimensional Hausdorff measure.  
\begin{proposition} \label{bound on H2}
Suppose a metric surface $U$ has the property that there are constants $C,R>0$ and a compact neighborhood $V$ such that $\sur(p,r)<Cr$,
for all $p\in V$ and all $r<R$.
Then, for $r<R$, any
$r$-ball in $V$ has $2$-dimensional
Hausdorff measure    greater than $C r^2$.
\end{proposition}

If the space is countably $2$-rectifiable, the Hausdorff  $2$-measure of an $R$-ball can be calculated by integrating
 from $0$ to $R$  the $1$-Hausdorff measure of the boundary of the $r$-ball in $dr$.  
 If the space is not countably $2$-rectifiable, the integral is always a lower bound (up to some factor), cf.  \cite{Fed69}.
Let $\H^k(X)$ be the $k$-dimensional Hausdorff measure of a metric space $X$.
We will make use of the following theorem. 
\begin{theorem}[{Federer, \cite[2.10.25]{Fed69}}]
Let $X$ be a metric space and let $f:X\to\R $ be a Lipschitz map. If $A\subset X$ and $k,m\geq 0$, then
$$ 
(\operatorname{Lip} f)^m\frac{\omega(k)\omega(m)}{\omega(k+m)}\H^{k+m}(A)
\geq \int_\R^*\H^k(A\cap f^{-1}\{r\})d\H^m(r)
,$$
where $ \int^*$ is the upper integral and $\omega(k)$ is the measure of the $k$-dimensional unit ball.
 \end{theorem}
\proof[Proof of Proposition \ref{bound on H2}] Using the theorem for $f(\cdot)=d(p,\cdot)$ (which is $1$-Lipschitz), $A=B(p,R)$, and $k=m=1$, we have
\begin{eqnarray*}
\frac{\omega(1)^2}{\omega(2)}\H^2(B(p,R))
 &\geq&\int_\R^*\H^1(B(p,R)\cap f^{-1}\{r\})d\H^1(r)\\
&=&\int_{[0,R]}^*\H^1(\partial B(p,t))dt
.\end{eqnarray*}
For the last equality, note that $f^{-1}\{r\}=\partial B(p,r)$.
Thus
\begin{equation}\label{accadue}
\H^2(B(p,R))\geq C_1 \int_{[0,R]}^*\H^1(\partial B(p,r))dr,\end{equation}
where $C_1$ is a suitable constant.

We claim that  $\H^1(\partial B(p,r))\geq Cr$. The rest of the subsection will be devoted to the demonstration of the claim. However, modulo this claim, the theorem is proved. Indeed, using it in  (\ref{accadue}) and integrating, we get
 $\H^2(B(p,R))\geq \dfrac{C}{2} R^2$. \qed
 
The reason behind the claim is that either $\partial B(p,r)$ has infinite length or it is a curve surrounding the ball $ B(p,r)$. If the measure is infinite there is nothing to prove. Consider the case when the measure is finite. Call  $\Sigma$ the exterior boundary of $ B(p,r)$, i.e., the boundary of the unbounded component of the complement of  $ B(p,r)$. Note that $\Sigma$ surrounds $ B(p,r)$, then if  $\Sigma$ were a curve, its $1$ dimensional Hausdorff measure would be its length. Thus the assertion of the claim follows from the bound on the surrounding function.

 To prove that $\Sigma:=\partial_{\text{ext}}B(p,r)$ is a curve, we want to use a general theorem \cite{Mazurkiewicz}:
\begin{theorem} [The Hahn-Mazurkiewicz theorem]
   A Hausdorff topological space is a continuous image of the unit
interval if and only if it is a Peano space, i.e., it is a compact, connected, locally connected
metric space.
 \end{theorem}
 To apply the theorem we only need to prove that $\Sigma$ is locally connected. 
By a corollary of the Phragm\'en-Brouwer  theorem, see \cite[page 106]{Whyburn},   since $\Sigma$  is a common boundary of two domains, it is a continuum.
In order to complete the proof of Proposition \ref{bound on H2} we just need to recall the following:
\begin{proposition} \label{continuum} Each  continuum  $\Sigma$ with $\H^1(\Sigma)<\infty$ is locally connected.\end{proposition}
A proof of the proposition can be argued using Theorem 12.1 in  \cite[page 18]{Whyburn}.
In what follows we give an alternative and easier proof.
\proof[Proof of Proposition \ref{continuum}]
Suppose that  $\Sigma$ is not locally connected. Hence there exist a point $p$ and a closed normal neighborhood $V$ of it such that any other neighborhood of $p$ contained in $ V$ is not connected.
 
\begin{lemma}\label{fatto1}
The closed set  $Z:=\cap\{S\;|\;p\in S, S\subset V, S {\rm \;clopen}\}$  is not a  neighborhood  of $p$. \end{lemma}
\proof
 Suppose $Z $ is a  neighborhood of $p$.
 Since $Z$ has to be disconnected,  there are $Z_1$ and $Z_2$ two closed (and therefore compact), disjoint  subsets of $Z$ such that $Z = Z_1 \cup Z_2$ and  $p\in Z_1$ but $p\notin Z_2$.

Since $V$ is normal, in $V$ there are  disjoint open neighborhoods $H_1$ and $H_2$ of $Z_1$ and $Z_2$ respectively. Let $H = H_1 \cup H_2$.

Since $  V \setminus H$ is a compact subset of $V\setminus Z$ there is a finite number of clopen subsets $K_1,\ldots,K_n$ of $V$ not containing  $p$ that cover $V \setminus H$.  Their union $K$ is also a clopen subset of $  V $, not containing  $p$ that covers   $V \setminus H$.
Clearly $K\cup H_2$ is a clopen subset of $V$ containing  $Z_2$ but not $p$.
\qed

Now fix a closed neighborhood $U\subset V$ of $p$ such that $c:=$dist$(U,\partial V)>0.$ 
By Lemma \ref{fatto1} there is a non-empty clopen set $Y$ of $V$ that  intersects $U$ but does not contain $p$. Since $\Sigma$ is connected and $U$ and $V$ are closed (and clearly different from $\Sigma$), $Y$ also intersects  $\partial U$ and $\partial V$ non-trivially.

\begin{lemma}\label{fatto2}
$\H^1(Y)\ge c$.\end{lemma}
\proof The function $\rho: Y \rightarrow \R$ defined by $\rho(y)=$dist$ (y,\partial V)$ is non-expanding.   Suppose there is a  point $\xi\in\R$  disconnecting  $\rho(Y)\subset\R$.
Then the set of all points of $Y$ with distance from $\partial V$ bigger than  $\xi$ is a clopen set of $Y$
 not intersecting $\partial V$.
Thus such a set is a proper clopen subset of $\Sigma$. This contradicts the fact that $\Sigma$ is connected.  Hence $\rho(Y)$ is a connected subset of the positive  real line, and moreover it contains $0$ and $c$. 
 Therefore the image of $\rho$ contains the interval $[0,c]$. Since $1$-Lipschitz maps do not increase Hausdorff measures and $\H([0,c])=c$,
we get $\H^1(Y)\ge c$.
 \qed
 
We can now conclude the proof of Proposition \ref{continuum}by contradicting the fact that 
$\H^1(\Sigma)<+ \infty$.
We will construct a sequence $Y_i$ of disjoint clopen subsets of $V$ with $\H^1(Y_i)\ge c$ for each $i$ and   arrive at a contradiction since  $\H^1(\Sigma) \ge \H^1(V) \ge \sum_i \H^1(Y_i)= + \infty$.


Put $U_1=U$, $V_1=V$ and $Y_1=Y$.
Inductively, consider $U_{j+1}:=U_j\setminus Y_j$ and $V_{j+1}:=V_j\setminus Y_j$, which are still closed. Using Lemma \ref{fatto1}, choose a clopen set $Y_ {j+1} $  of $V_{j+1}$ that does not contain $p$ but meets  $ U_{j+1}$, hence it meets also  $\partial U_{j+1}$ and  $\partial V_{j+1}$.

Note that  since  $V_j$ is a clopen subset of $V$ then  $V_j \setminus \partial V$ is open and so $\partial V_j \subset \partial V$. Similarly $ \partial U_j \subset \partial U $ and hence  dist$(\partial U_j,\partial V_j)\ge$dist$(\partial U,\partial V)\ge c$.

As  for $Y=Y_1$ we have that  $\H^1(Y_j)\ge c$.   
\qed

 \bibliography{bilip_homog_geod_surf_bib}
\bibliographystyle{amsalpha}
 
 
\vskip 1in

\parbox{3.5in}{Enrico Le Donne:\\
~\\
Department of Mathematics\\
Yale University\\
New Haven, CT 06520\\
enrico.ledonne@yale.edu}

\end{document}